\documentclass{article}
\catcode`\@=11

\newtheorem{proposition}{Proposition}

\newtheorem{example}{Example}
\newtheorem{remark}{Remark}

\newtheorem{corollary}{Corollary}
\def\demo{\noindent{\bf Proof .-}}
\def\section{\@startsection {section}{1}{\z@}{-3.5ex plus -1ex
minus-.2ex}{2.3ex plus .2ex}{\normalsize\bf}}

\def\gcd{{\rm gcd}\,}

\begin{document}
\begin{center}
{\Large\bf \textsc{Determinantal ideals and monomial curves in the three-dimensional space}}
\end{center}
\vskip.5truecm
\begin{center}
{Margherita Barile\footnote{Partially supported by the Italian Ministry of Education, University and Research.}\\ Dipartimento di Matematica, Universit\`{a} di Bari,Via E. Orabona 4,\\70125 Bari, Italy}
\end{center}
\vskip1truecm
\noindent
{\bf Abstract} We show that the defining ideal of every monomial curve in the affine or projective three-dimensional space can be set-theoretically defined by three binomial equations, two of which set-theoretically define a determinantal ideal generated by the 2-minors of a $2\times 3$ matrix with monomial entries. 

\section*{Introduction} Let $K$ be an algebraically closed field, and let $R$ be a polynomial ring in $n$ indeterminates over $K$. If $I$ is a proper ideal of $R$, then $I$ is called a {\it set-theoretic complete intersection} if the corresponding variety $V(I)$ in the affine space $K^n$ can be defined by a system containing the least possible number $s$ of equations: in this case $s=$\, height\,$(I)$. More generally, if we can take $s\leq$\,height\,$(I)+1$, then $I$  is called an {\it almost set-theoretic complete intersection}. We recall that an ideal is a(n ideal-theoretic) {\it complete intersection} if it is generated by height\,$(I)$ elements. \newline
The set-thoretic complete intersection property has been intensively studied for ideals generated by binomials, such as the determinantal ideals (see, e.g., \cite{Ba}, \cite{BS} and \cite{BV}), and the ideals defining toric varieties (see, e.g., \cite{BMT}), which include the special class of monomial curves (see, e.g., \cite{B}, \cite{H}, \cite{M}, \cite{P}, and \cite{T}). Monomial curves in $K^3$ are known to always be set-thoretic complete intersections (see \cite{B} and \cite{H}), and the same is true for all monomial curves in ${\bf P}^3$ (see \cite{M}) if char\,$K>0$: the extension to the characteristic zero case is a long-standing open problem, besides some special cases (see, e.g., \cite{RV}). What makes the problem so hard is the fact that, in characteristic zero, a monomial curve (and, more generally, a toric variety) is  not a set-theoretic complete intersection on ``simple"  (i.e., on {\it binomial}) equations, except in the trivial case where it is a complete intersection (see \cite{BMT}, Theorem 4).  Even for determinantal varieties, the minimum number of defining equations is in general provided by non-binomial ones (see \cite{B}, and  also \cite{BS}, Theorem 2 together with \cite{BV}, Section 5.E or with \cite{Br}, Section 2): this is also true for monomial curves in $K^3$ and for arithmetically Cohen-Macaulay curves in ${\bf P}^3$, which are defined by the vanishing of the 2-minors of a $2\times 3$ matrix (see \cite{B} or \cite{H}, and \cite{BSV} or \cite{RV}). \newline
In this paper we first describe a class of ideals, generated by the maximal minors of a two-row matrix, whose variety is defined by the vanishing of a proper subset of the generating minors; then we show that the defining ideal of any monomial curve in $K^3$ or ${\bf P}^3$ contains the ideal of 2-minors of a $2\times 3$ matrix that is a set-theoretic complete intersection on two of these minors. Moreover, these two minors and an additional binomial define the curve set-theoretically. This is an attempt to give a unifying feature of affine and projective monomial curves, and could help us to shed some  light on the conjecture according to which every projective curve in ${\bf P}^3$ is a set-theoretic complete intersection.  This also provides a new proof to the following well-known result: the defining ideal of every monomial curve in the three-dimensional space contains a complete intersection ideal of height 2 generated by binomials. The problem of finding complete intersection ideals in toric ideals has been recently treated in \cite{CCD}. 
\section{A result on determinantal ideals}
Let $A=(\alpha_{ij})$ be a $2\times r$ matrix with entries in $R$, where $r\geq 3$. Moreover, let $J$ be the ideal of $R$ generated by the 2-minors of $A$. For all distinct indices $i,j$ let $\Delta_{ij}$ denote the minor of $A$ formed by the $i$-th and the $j$-th column. Finally, for all $k=1,\dots, r$,  let $J_k$ be the ideal of $R$ generated by the set $\{\Delta_{ik}\vert 1\leq i\leq r,\ i\neq k\}$; in other words, $J_k$ is generated by all 2-minors of $A$ which involve the $k$-th column.
\begin{proposition}\label{ik} $\sqrt J=\sqrt{J_k}$ if and only if $J\subset\sqrt{(\alpha_{1k}, \alpha_{2k})}$.
\end{proposition}
\demo Since $J_k\subset (\alpha_{1k}, \alpha_{2k})$, the {\it  only if} part is trivial. We prove the {\it if} part. For the sake of simplicity, and without loss of generality, we assume that $k=1$. Thus our assumption is 
\begin{equation}\label{0}J\subset\sqrt{(\alpha_{11}, \alpha_{21})},
\end{equation}
\noindent
 and from this we want to deduce that $\sqrt J=\sqrt{J_1}$, or,  equivalently, that $V(J)=V(J_1)$. Since $J_1\subset J$, it suffices to show that $V(J)\supset V(J_1)$. Let ${\bf v}\in K^n$ be a point where all elements of $J_1$ vanish, i.e., such that 
\begin{equation}\label{1} \alpha_{11}({\bf v})\alpha_{2j}({\bf v})-\alpha_{1j}({\bf v})\alpha_{21}({\bf v})=\Delta_{1j}({\bf v})=0\qquad\mbox{for all }j=2,\dots, r.
\end{equation}
\noindent
Let
$$B=\left(\begin{array}{ccccc}
&&A&\\
\hline\\
\alpha_{11}&\alpha_{12}&\cdots&\cdots&\alpha_{1r}
\end{array}
\right).
$$
\noindent
Since the first and the last row of $B$ are equal, all 3-minors of $B$ vanish. In particular, for all indices $i,j$ such that $1<i<j\leq r$ we have
$$0=\left\vert
\begin{array}{ccc}
\alpha_{11}&\alpha_{1i}&\alpha_{1j}\\
\alpha_{21}&\alpha_{2i}&\alpha_{2j}\\
\alpha_{11}&\alpha_{1i}&\alpha_{1j}
\end{array}
\right\vert
=\alpha_{11}\Delta_{ij}-\alpha_{1i}\Delta_{1j}+\alpha_{1j}\Delta_{1i},$$
\noindent
whence
$$\alpha_{11}\Delta_{ij}\in J_1,$$
\noindent
and, consequently, due to the choice of ${\bf v}$, 
$$\alpha_{11}({\bf v})\Delta_{ij}({\bf v})=0.$$
\noindent
We have to show that $\Delta_{ij}({\bf v})=0$ for all indices $i,j$. This is certainly true if $\alpha_{11}({\bf v})\ne0$. So assume that $\alpha_{11}({\bf v})=0$. Then for all $j=2,\dots, r$, (\ref{1}) implies that 
$\alpha_{1j}({\bf v})\alpha_{21}({\bf v})=0$. If $\alpha_{21}({\bf v})=0$, then (\ref{0}) implies that 
\begin{equation}\label{2}
\Delta_{ij}({\bf v})=0\qquad\mbox{for all }i,j. 
\end{equation}
\noindent
Otherwise $\alpha_{1j}({\bf v})=0$ for all $j=2,\dots, r$, so that the whole first row of $A$ vanishes at ${\bf v}$, and (\ref{2}) holds in this case, too. This proves our claim.
\section{Matrices with monomial entries}
In this section we shall apply Proposition \ref{ik} to the case where the entries of $A$ are {\it monomials}, i.e., products of nonnegative powers of indeterminates. A difference of two distinct monomials will be called a {\it binomial}. We shall say that a matrix with monomial entries is {\it simple} if none of its 2-minors has a non constant monomial factor.  
\begin{example}{\rm In the polynomial ring $R=K[a,b,c,d]$ consider the matrix
$$A=\left(\begin{array}{ccc}
a^md^n& b^p&c^q\\
b^r&a^s&d^t
\end{array}\right),$$
\noindent
where $m,n,p,q,r,s,t$ are nonnegative integers. 
The ideal generated by the 2-minors of $A$ is 
$$J=(a^{m+s}d^n-b^{p+r},\ a^md^{n+t}-b^rc^q,\ b^pd^t-a^sc^q),$$
\noindent
so that 
$$J\subset(a,b)\subset\sqrt{(a^s, b^p)}.$$
\noindent
By Proposition \ref{ik} it follows that
$$\sqrt{J}=\sqrt{J_2}=\sqrt{(a^{m+s}d^n-b^{p+r}, b^pd^t-a^sc^q )}.$$
\noindent
}
\end{example}
\begin{remark}{\rm Let $J$ be the ideal of 2-minors of a $2\times3$ matrix. Then $J$ is generated by three elements of $R$, and,
according to Eagon and Northcott (\cite{EN}, Theorem 3), height\,$J\leq2$. Hence, if equality holds, $J$ is an almost set-theoretic complete intersection. In general, it is not a set-theoretic complete intersection; for instance, it is not if the entries of the matrix are pairwise distinct indeterminates. In this case height\,$J=2$, but $J$ cannot be generated, up to radical, by less than three polynomials (see \cite{BS}, Theorem 2). This shows the interest of
Proposition \ref{ik}, which describes a class of $2\times 3$ matrices with entries in $R$ such that the ideal generated by its 2-minors  is generated by two polynomials up to radical. A different class was characterized by Robbiano and Valla, who showed the following result.
}\end{remark}
\begin{proposition}\label{valla}{\rm(\cite{RV}, Theorem 2.2)}. Let $J$ be the ideal of $R$ generated by the $2$-minors of the matrix 
$$A=\left(\begin{array}{ccc}
x&y^m&z\\
y^n&s&t
\end{array}\right),
$$
\noindent
where $m,n$ are nonnegative integers. Then there are $f,g\in R$ such that $\sqrt J=\sqrt{(f,g)}$.
\end{proposition}
\par\medskip\noindent
Propositions \ref{ik} and \ref{valla} can be summed up as follows.
\begin{corollary}\label{corollary} Let $A$ be a simple $2\times3$ matrix whose entries are monomials of $R=K[a,b,c,d]$.  Let $J$ be the ideal generated by its $2$-minors. Then there are $f,g\in R$ such that $\sqrt{J}=\sqrt{(f,g)}$.
\end{corollary} 
\demo It suffices to show that $A$ fulfils the assumptions of Proposition \ref{ik} or Proposition \ref{valla}. The assumption of Proposition \ref{ik} is certainly true if one of the entries of $A$ is equal to 1. So assume that this is not the case. Then note that, according to the simplicity condition, each entry is coprime with respect to the entries lying in the same row. Hence, up to permuting  columns and renaming the indeterminates, the first row of matrix $A$ is of the following form
$$\left(\begin{array}{ccc}
a^m&b^n&c^pd^q
\end{array}\right),
$$
\noindent
 where $m,n$ are positive integers and $p, q$ are nonnegative integers, not both zero. The entries of the second row are, of course, of the same form. Moreover, the simplicity condition implies that every entry is coprime with respect to the entries lying in the same column. Hence the possible forms of $A$ are the ones in the list below (up to interchanging $c$ and $d$). We assume that each monomial contains at least one positive exponent. 
\begin{list}{}{}
\item{(i)} 
 $$\left(\begin{array}{ccc}
a^m&b^n&c^pd^q\\
b^r&c^sd^t&a^u
\end{array}\right)
$$
\noindent
\item{(ii)} 
 $$\left(\begin{array}{ccc}
a^m&b^n&c^p\\
b^r&a^sc^t&d^u
\end{array}\right)
$$
\noindent
\item{(iii)} 
 $$\left(\begin{array}{ccc}
a^m&b^n&c^p\\
b^r&c^s&a^td^u
\end{array}\right)
$$
\noindent
\item{(iv)} 
 $$\left(\begin{array}{ccc}
a^m&b^n&c^pd^q\\
b^rc^s&d^t&a^u
\end{array}\right)$$
\noindent
\item{(v)} 
 $$\left(\begin{array}{ccc}
a^m&b^n&c^p\\
b^rc^s&a^t&d^u
\end{array}\right)
$$
\noindent
\item{(vi)}
$$\left(\begin{array}{ccc}
a^m&b^n&c^pd^q\\
c^r&a^sd^t&b^u
\end{array}\right)
$$
\noindent
\item{(vii)}
$$\left(\begin{array}{ccc}
a^m&b^n&c^pd^q\\
c^r&d^s&a^tb^u
\end{array}\right)
$$
\noindent
\item{(viii)} 
 $$\left(\begin{array}{ccc}
a^m&b^n&c^p\\
c^r&a^s&b^td^u
\end{array}\right)
$$
\noindent
\item{(ix)}
$$\left(\begin{array}{ccc}
a^m&b^n&c^pd^q\\
c^rd^s&a^u&b^t
\end{array}\right)
$$
 \end{list}
\noindent
Proposition \ref{valla} can be applied to (i), (ii), (iii), (iv), (v), (vi), (viii) and (ix). It can be applied in case (vii) if one of the exponents $p,q,t,u$ is zero. Otherwise we have $J\subset (a,c)$ (which is true if $t,p>0$) and $J\subset (b,d)$ (which is true if $q,u>0$), so that Proposition \ref{ik} applies. This completes the proof. 
\par\medskip\noindent
\begin{remark}{\rm Note that (ii), (iv), (v) and (vi)  may also fulfil the assumption of Proposition \ref{ik}. This is true in case (ii) if $s>0$, since then $J\subset(a,b)$; in case (iv) if $r=0$ and $p>0$, or $q=0$ and $s>0$, since then $J\subset (a,c)$; in case (v) if $r>0$, since then $J\subset (a,b)$; in case (vi) if $s=0$ and $q>0$, or $p=0$ and $t>0$, since then $J\subset (b,d)$.    In these special cases, if we are asked to explicitly exhibit the polynomials $f,g$ generating $J$ up to radical, Proposition \ref{ik} is in fact more convenient than Proposition \ref{valla}: the former yields two binomials, the latter a binomial and polynomial of a more complicated form. We compare the two approaches in the next examples.}
\end{remark} 
\begin{example}\label{example2}{\rm Let 
$$A=\left(\begin{array}{ccc}
a^m&b^n&c^p\\
b^r&c^s&a^u
\end{array}\right),$$
\noindent
with positive exponents, which is a matrix of type (i) for $t=q=0$. Then $J=(a^mc^s-b^{n+r},\ a^{m+u}-b^rc^p,\ a^ub^n-c^{p+s})$. Two polynomials $f,g\in K[a,b,c]$ generating  $J$ up to radical are given by Valla through an explicit formula in \cite{V}, Section 3:
\begin{eqnarray*} f&=&a^mc^s-b^{n+r}\\
g&=&\sum_{k=0}^{n+r}(-1)^{n+r-k}{{n+r}\choose{k}}a^{ku+\tau_k m}b^{\sigma_k}c^{(n+r-k)(p+s)+\tau_ks-ns},
\end{eqnarray*} 
\noindent
where $\tau_k$ and $\sigma_k$ are the quotient and the remainder of the Euclidean division of $kn$ by $n+r$. Here Proposition \ref{ik} cannot be applied. It only works in the ``degenerate" case, when one of the exponents is set equal to zero. If, e.g., $n=0$, then  $f=a^mc^s-b^r$ and we can take $g=a^u-c^{p+s}$. In this special case one can easily check that we even have $J=(f,g)$.   
} 
\end{example}
\begin{example}\label{example3}{\rm Let  
 $$A=\left(\begin{array}{ccc}
a^{2u-1}&b^n&c^p\\
b^rc^s&d^t&a^u
\end{array}\right),
$$
\noindent
so that $J=(a^{2u-1}d^t-b^{n+r}c^s,\ a^ub^n-c^pd^t,\ a^{3u-1}-b^rc^{p+s})$. Suppose that the exponents are all positive.
Applying the general costructive method described in \cite{RV} requires more than one step for all $u>1$, and one of the equations we obtain is non-binomial. Proposition \ref{ik} immediately gives  $\sqrt J=\sqrt{(a^ub^n-c^pd^t,\ a^{3u-1}-b^rc^{p+s})}$. 
}
\end{example}
\section{Complete intersections and monomial curves}
There is a well-known connection between the defining ideals of affine and projective monomial curves and determinantal ideals of simple matrices. Consider the affine curve of $K^3$ parametrized by
$$x_1=\xi^{\alpha},\ x_2=\xi^{\beta},\ x_3=\xi^{\gamma}.$$
\noindent
Its defining ideal in $K[x_1, x_2, x_3]$ is the ideal of 2-minors of a matrix
$$A=\left(\begin{array}{ccc}
x_1^m&x_2^n&x_3^p\\
x_2^r&x_3^s&x_1^t
\end{array}\right),$$
\noindent
which has the form given in Example \ref{example2}. By Proposition \ref{valla} it is a set-theoretic complete intersection; see also \cite{B} and \cite{H}. This property does not extend to the projective case. Consider the projective curve of ${\bf P}^3$ parametrized by 
\begin{equation}\label{parameter} C: x_0=\xi^{\delta},\ x_1=\xi^{\epsilon_1}\omega^{\delta-\epsilon_1},\ x_2=\xi^{\epsilon_2}\omega^{\delta-\epsilon_2},\ x_3=\omega^{\delta},\end{equation}
\noindent
where $\delta,\epsilon_1, \epsilon_2$ are positive integers. The defining ideal $I=I(C)$ of $C$ in $R=K[x_0, x_1, x_2,x_3]$ is not necessarily the determinantal ideal of a $2\times 3$ matrix. In fact it is generated by three elements if and only if it is arithmetically Cohen-Macaulay, in which case it is generated by the 2-minors of a matrix
$$A=\left(\begin{array}{ccc}
x_1^m&x_2^n&x_0^px_3^q\\
x_2^r&x_0^sx_3^t&x_1^u
\end{array}\right),$$
\noindent
which is of type (i). This result was found by Bresinsky, Schenzel and Vogel in \cite{BSV}, see also \cite{RV}, Theorem 2.1.  According to Proposition \ref{valla}, it follows that $I$ is a set-theoretic complete intersection. It was proven by Moh \cite{M} that the set-theoretic complete intersection property extends to {\it all} projective monomial curves $C$ if the characteristic of $K$ is positive, whereas the question is open in characteristic zero. Our aim is to show a general property that links all ideals $I$ to determinantal ideals which are set-theoretic complete intersections. A binomial will be called {\it monic} if one of its two monomials is a power of an indeterminate. Before stating the main result, we need to introduce some notation, which is referred to the projective curve $C$ defined by (\ref{parameter}). Set $\phi_i=\delta-\epsilon_i$ for $i=1,2$. Up to exchanging the parameters $\xi$ and $\omega$ we may assume that $\epsilon_1\geq\epsilon_2$ (or, equivalently, $\phi_2\geq\phi_1$). Then let
\begin{eqnarray*}
\partial&=&\gcd(\delta, \epsilon_1)=\gcd(\delta, \phi_1)\\
\varepsilon&=&\gcd(\epsilon_1, \epsilon_2)\\
\varphi&=&\gcd(\phi_1, \phi_2).
\end{eqnarray*}
\noindent
Further set
$$\begin{tabular}{ccc}
$\delta^{\ast}=\displaystyle\frac\delta\partial;$&$\epsilon_1^{\ast}=\displaystyle\frac{\epsilon_1}\partial;$&$\phi_1^{\ast}=\displaystyle\frac{\phi_1}\partial$;\\\\
$m=\displaystyle\frac{\phi_2}{\varphi};$&$p=\displaystyle\frac{\phi_1}{\varphi};$&$n=m-p$;\\\\
$u=\displaystyle\frac{\epsilon_1}{\varepsilon};$&$v=\displaystyle\frac{\epsilon_2}{\varepsilon};$&$w=u-v$.\\
\end{tabular}$$
\noindent
Note that $n$ and $w$ are nonnegative integers. Consider the following three binomials
\begin{eqnarray}
\label{f}f&=&x_1^{\delta^{\ast}}-x_0^{\epsilon_1^{\ast}}x_3^{\phi_1^{\ast}};\\
\label{f1}f_1&=&x_1^m-x_0^nx_2^p;\\
\label{f2}f_2&=&x_2^u-x_1^vx_3^w.
\end{eqnarray}
\noindent
Since
\begin{eqnarray}
\label{5}\epsilon_1\delta^{\ast}&=&\epsilon_1^{\ast}\delta\quad\quad\quad(\mbox{equivalently: }\phi_1\delta^{\ast}=\phi_1^{\ast}\delta);\\
\label{6}\epsilon_1m&=&\delta n+\epsilon_2p\quad(\mbox{equivalently: }\phi_1m=\phi_2p);\\
\label{7}\epsilon_2u&=&\epsilon_1v\quad\quad\quad(\mbox{equivalently: }\phi_2u=\phi_1v+\delta w),
\end{eqnarray}
\noindent
we have that $f, f_1, f_2\in I(C)$. 
\begin{proposition}\label{ci} The monomial curve $C$ in ${\bf P}^3$ is set-theoretically defined by the binomials $f,f_1, f_2$.\end{proposition} 
\demo  It suffices to prove that $V(f,f_1,f_2)\subset C$. 
Let ${\bf x}=(x_0,x_1,x_2,x_3)\in{\bf P}^3$ be a common zero of $f, f_1, f_2$. Set $x_0=\xi$ and $x_3=\omega$. We will show that, after a suitable change of parameters, ${\bf x}$ fulfils parametrization (\ref{parameter}). This is true if $x_0=0$ or $x_3=0$. If $x_0=0$, then $f({\bf x})=0$ implies that $x_1=0$, and then $f_2({\bf x})=0$ implies that $x_2=0$; we thus can take $\xi=0$ and set $\omega$ equal to any $\delta$-th root of $x_3$. Similarly, if $x_3=0$, we conclude  that $x_1=x_2=0$, so that we can take $\omega=0$ and set $\xi$ equal to any $\delta$-th root of $x_0$. So assume that $x_0\ne 0$ or $x_3\ne 0$.
Since $f({\bf x})=0$, we have that
\begin{equation}\label{7'}x_1^{\delta^{\ast}}=x_0^{\epsilon_1^{\ast}}x_3^{\phi_1^{\ast}}=\xi^{\epsilon_1^{\ast}}\omega^{\phi_1^{\ast}}.\end{equation}
\noindent 
Let $\bar\xi, \bar\omega\in K$ be such that ${\bar\xi}^{\delta}=\xi$, and ${\bar\omega}^{\delta}=\omega$, i.e.,  
\begin{equation}\label{8}
x_0={\bar\xi}^{\delta},\quad\mbox{and}\quad x_3={\bar\omega}^{\delta}.
\end{equation}
\noindent
Then, by (\ref{7'}) and (\ref{5}),
$$x_1^{\delta^{\ast}}={\bar\xi}^{\delta\epsilon_1^{\ast}}{\bar\omega}^{\delta\phi_1^{\ast}}=
{\bar\xi}^{\delta^{\ast}\epsilon_1}{\bar\omega}^{\delta^{\ast}\phi_1},$$
\noindent
whence
\begin{equation}\label{9}
x_1={\bar\xi}^{\epsilon_1}{\bar\omega}^{\phi_1}r,
\end{equation}
\noindent
for some $r\in K$ such that
\begin{equation}\label{10}
r^{\delta^{\ast}}=1.
\end{equation}
\noindent
From $f_2({\bf x})=0$ it follows that
$$
x_2^u=x_1^vx_3^w,
$$
\noindent
so that, in view of (\ref{9}) and (\ref{8}),
\begin{eqnarray*}
x_2^u&=&{\bar\xi}^{\epsilon_1v}{\bar\omega}^{\phi_1v}r^v{\bar\omega}^{\delta w}=
{\bar\xi}^{\epsilon_1v}{\bar\omega}^{\phi_1v+\delta w}r^v\\
&=&{\bar\xi}^{\epsilon_2u}{\bar\omega}^{\phi_2u}r^v,
\end{eqnarray*}
\noindent
where the last equality follows from (\ref{7}). Hence
\begin{equation}\label{11}
x_2={\bar\xi}^{\epsilon_2}{\bar\omega}^{\phi_2}s,
\end{equation}
\noindent
for some $s\in K$ such that
\begin{equation}\label{12}
s^u=r^v.
\end{equation}
\noindent
From $f_1({\bf x})=0$ it follows that
$$
x_1^m=x_0^nx_2^p.
$$
\noindent
By (\ref{9}), (\ref{8}) and (\ref{11}) this is equivalent to
$${\bar\xi}^{\epsilon_1m}{\bar\omega}^{\phi_1m}r^m={\bar\xi}^{\delta n+\epsilon_2p}{\bar\omega}^{\phi_2p}s^p.$$
\noindent
In view of (\ref{6}), cancelling equal powers of $\bar\xi$ and $\bar\omega$ yields
\begin{equation}\label{13}
r^m=s^p.
\end{equation}
From the definition of $m,p,u,v$ we have
\begin{eqnarray*}
m\varphi+v\varepsilon&=&\phi_2+\epsilon_2=\delta\\
p\varphi+u\varepsilon&=&\epsilon_1+\phi_1=\delta.
\end{eqnarray*}
\noindent
Hence, by (\ref{12}) and (\ref{13}), it holds that
\begin{equation}\label{14}
s^{\delta}=s^{p\varphi+u\varepsilon}=r^{m\varphi+v\varepsilon}=r^{\delta}=1,
\end{equation}
\noindent
where the last equality is a consequence of (\ref{10}). Thus both $s$ and $r$ are $\delta$-th roots of unity. Let $\eta$ be a primitive $\delta$-th root of unity. Then
\begin{equation}\label{15} r=\eta^e,\quad\mbox{and}\quad s=\eta^f
\end{equation}
for suitable integers $e$, $f$. Since, in view of (\ref{10}), the order of $\eta^e$ 
as a root of unity is a divisor of $\delta^{\ast}$, it follows that $\partial$ divides $e$. Moreover, replacing (\ref{15}) in (\ref{12}) yields
$$\eta^{uf}=\eta^{ve},$$
\noindent
whence
$$\delta\mbox{ divides } uf-ve=\frac{\epsilon_1f-\epsilon_2e}{\varepsilon},$$
\noindent
so that, in particular,
\begin{equation}\label{16}
\epsilon_1f\equiv\epsilon_2e\qquad\qquad(\mbox{mod}\ \delta)
\end{equation}
\noindent
Let $a=\displaystyle\frac{e}{\partial}\bar\epsilon_1$, where $\bar\epsilon_1$ is an integer such that
\begin{equation}\label{17}
\bar\epsilon_1\epsilon_1^{\ast}\equiv1\qquad\qquad(\mbox{mod}\ \delta)
\end{equation}
\noindent
Such $\bar\epsilon_1$ exists because $\gcd(\epsilon_1^{\ast},\delta)=1$. Then
\begin{equation}\label{18}
\eta^{a\epsilon_1}=\eta^{e\bar\epsilon_1\epsilon_1/\partial}=\eta^{e\bar\epsilon_1\epsilon_1^{\ast}}=\eta^e=r,
\end{equation}
\noindent
where the last two equalities follow from (\ref{17}) and (\ref{15}) respectively. Moreover
\begin{equation}\label{19}
\eta^{a\epsilon_2}=\eta^{e\bar\epsilon_1\epsilon_2/\partial}=\eta^{f\bar\epsilon_1\epsilon_1/\partial}=\eta^{f\bar\epsilon_1\epsilon_1^{\ast}}=\eta^f=s,
\end{equation}
\noindent
where the second, the fourth and the fifth equality follow from (\ref{16}),  (\ref{17}) and (\ref{15}) respectively. Set 
\begin{equation}\label{20}
\tilde\xi=\bar\xi\eta^a.
\end{equation}
\noindent
We then replace (\ref{20}) in (\ref{8}), (\ref{9}) and (\ref{11}). We obtain
\begin{eqnarray}
\label{21}x_0&=&\bar\xi^{\delta}=\bar\xi^{\delta}\eta^{a\delta}=\tilde\xi^{\delta},\\
\label{22}\mbox{by (\ref{18})}\quad\quad x_1&=&\bar\xi^{\epsilon_1}\bar\omega^{\phi_1}r=\bar\xi^{\epsilon_1}\bar\omega^{\phi_1}\eta^{a\epsilon_1}=\tilde\xi^{\epsilon_1}\bar\omega^{\phi_1},\\
\label{23}\mbox{ and, by (\ref{19})}\quad\quad x_2&=&\bar\xi^{\epsilon_2}\bar\omega^{\phi_2}s=\bar\xi^{\epsilon_2}\bar\omega^{\phi_2}\eta^{a\epsilon_2}=\tilde\xi^{\epsilon_2}\bar\omega^{\phi_2}.
\end{eqnarray}
\noindent
Equalities (\ref{21})--(\ref{23}), together with (\ref{8}), show that $\tilde\xi$ and $\bar\omega$ are the parameters which allow us to represent ${\bf x}$ in the form (\ref{parameter}). This completes the proof. \par\medskip\noindent
Now consider the following matrix with monomial entries:
$$A=\left(\begin{array}{ccc}
x_1^{\min(\delta^{\ast}, m)}&x_0^{\max(\epsilon_1^{\ast}-n,0)}x_3^{\phi_1^{\ast}}&x_0^{\max(n-\epsilon_1^{\ast},0)}x_2^p\\
x_0^{\min(\epsilon_1^{\ast}, n)}&x_1^{\max(\delta^{\ast}-m,0)}&x_1^{\max(m-\delta^{\ast},0)}
\end{array}\right).$$
\noindent
One can easily check that $A$ is simple. 
Let $J$ be the ideal generated by the 2-minors of $A$. In the sequel we will throughout refer to the projective curve $C$ given in (\ref{parameter}). 
\begin{corollary}\label{corollary2}
We have that $C=V(M_1, M_2, f_2)$, where  $M_1$, $M_2$ are two minors of $A$ such that $V(M_1,M_2)=V(J)$. 
\end{corollary} 
\demo
  We have to distinguish between different cases. In each case we will show that $A$ fulfils the assumption of Proposition \ref{ik} with respect to the second column ($k=2$) or the third column ($k=3$). We will choose $M_1$, $M_2$ accordingly,  so as to have $V(M_1, M_2)=V(J)$.
\begin{list}{}{}
\item{Case (I):} $\epsilon_1^{\ast}>n$, $\delta^{\ast}>m$. Then 
$$A=\left(\begin{array}{ccc}
x_1^m&x_0^{\epsilon_1^{\ast}-n}x_3^{\phi_1^{\ast}}&x_2^p\\
x_0^n&x_1^{\delta^{\ast}-m}&1
\end{array}\right).$$
\noindent In this case we can take $M_1=\Delta_{13}$ and $M_2=\Delta_{23}$.
\item{Case (II):} $\epsilon_1^{\ast}>n$, $\delta^{\ast}\leq m$. Then 
$$A=\left(\begin{array}{ccc}
x_1^{\delta^{\ast}}&x_0^{\epsilon_1^{\ast}-n}x_3^{\phi_1^{\ast}}&x_2^p\\
x_0^n&1&x_1^{m-\delta^{\ast}}
\end{array}\right).$$
\noindent In this case we can take $M_1=\Delta_{12}$ and $M_2=\Delta_{23}$.
\item{Case (III):} $\epsilon_1^{\ast}\leq n$, $\delta^{\ast}>m$. Then 
$$A=\left(\begin{array}{ccc}
x_1^m&x_3^{\phi_1^{\ast}}&x_0^{n-\epsilon_1^{\ast}}x_2^p\\
x_0^{\epsilon_1^{\ast}}&x_1^{\delta^{\ast}-m}&1
\end{array}\right).$$
\noindent In this case we can take $M_1=\Delta_{13}$ and $M_2=\Delta_{23}$.
\item{Case (IV):} $\epsilon_1^{\ast}\leq n$, $\delta^{\ast}\leq m$. Then 
$$A=\left(\begin{array}{ccc}
x_1^{\delta^{\ast}}&x_3^{\phi_1^{\ast}}&x_0^{n-\epsilon_1^{\ast}}x_2^p\\
x_0^{\epsilon_1^{\ast}}&1&x_1^{m-\delta^{\ast}}
\end{array}\right).$$
\noindent In this case we can take $M_1=\Delta_{12}$ and $M_2=\Delta_{23}$.
\end{list}
Now, in each of the above cases,
\begin{eqnarray}\label{d12}\Delta_{12}&=&f\\
\label{d13}\Delta_{13}&=&f_1,
\end{eqnarray}
\noindent
 whence $f, f_1\in J$, so that 
\begin{equation}\label{vf} V(J)=V(M_1,M_2)\subset V(f_1, f_2).
\end{equation}
\noindent
Moreover:
\begin{list}{}{}
\item{} in case (I)
$$\Delta_{23}=x_0^{\epsilon_1^{\ast}-n}x_3^{\phi_1^{\ast}}-x_1^{\delta^{\ast}-m}x_2^p;
$$
\item{} in case (II)
$$\Delta_{23}=x_1^{m-\delta^{\ast}}x_0^{\epsilon_1^{\ast}-n}x_3^{\phi_1^{\ast}}-x_2^p;
$$
\item{} in case (III)
$$\Delta_{23}=x_3^{\phi_1^{\ast}}-x_0^{n-\epsilon_1^{\ast}}x_1^{\delta^{\ast}-m}x_2^p;
$$
\item{} in case (IV)
$$\Delta_{23}=x_3^{\phi_1^{\ast}}x_1^{m-\delta^{\ast}}-x_0^{n-\epsilon_1^{\ast}}x_2^p.
$$
\end{list}
\noindent
From (\ref{5}) and (\ref{6}) it follows that 
$$\epsilon_1(\delta^{\ast}-m)=(\epsilon_1^{\ast}-n)\delta-\epsilon_2p,$$
which implies that in all cases $\Delta_{23}\in I(C)$. Thus $M_2\in I(C)$, and, in view of (\ref{d12}) and (\ref{d13}), $M_1\in\{f,f_1\}\subset I(C)$. Hence
$$C\subset V(M_1, M_2, f_2)\subset V(f, f_1, f_2)=C,$$
\noindent
where the second and the third inclusion follow from (\ref{vf}) and Proposition \ref{ci} respectively. 

This completes the proof.
 \par\medskip\noindent
\begin{example}\label{example4}{\rm   Consider the projective curve
$$C: C: x_0=\xi^4,\ x_1=\xi^3\omega,\ x_2=\xi\omega^3,\ x_3=\omega^4,$$
which is known as the {\it rational quartic} and is conjectured to be a set-theoretic complete intersection. Its defining ideal $I(C)$ is minimally generated by $x_1^2x_3-x_0x_2^2, x_1^3-x_0^2x_2, x_2^3-x_1x_3^2, x_0x_3-x_1x_2$. Curve $C$ fulfils case (I) in the proof 
of Corollary \ref{corollary2}, and corresponds to the matrix
$$A=\left(\begin{array}{ccc}
x_1^3&x_0x_3&x_2\\&&\\
x_0^2&x_1&1
\end{array}\right).$$
\noindent
Therefore  curve $C$ is set-theoretically defined  by
\begin{eqnarray*}
\Delta_{13}&=&x_1^3-x_0^2x_2\\
\Delta_{23}&=&x_0x_3-x_1x_2\\
f_2&=&x_2^3-x_1x_3^2.
\end{eqnarray*}
\noindent
Note that matrix $A$ satisfies the assumption of Proposition \ref{ik} with respect to the first and the third column. Hence both the pairs of minors $\Delta_{12}$, $\Delta_{13}$ and $\Delta_{13}$, $\Delta_{23}$ set-theoretically define the corresponding determinantal variety $V(J)$, $J$ being the ideal generated by the 2-minors of $A$.  It follows that $C$ is also set-theoretically defined by
\begin{eqnarray*}
f=\Delta_{12}&=&x_1^4-x_0^3x_3\\
f_1=\Delta_{13}&=&x_1^3-x_0^2x_2\\
f_2&=&x_2^3-x_1x_3^2.
\end{eqnarray*}
This is a well-known result on the rational quartic, which is also recalled in \cite{RV}, p.~391.
}\end{example}
\begin{remark}{\rm  The fact that every monomial curve in the three-dimensional space is set-theoretically defined by three monic binomials had already been proved in \cite{BM}. In that paper, however, the proof was not based on closed formulas like (\ref{f}), (\ref{f1}) and (\ref{f2}), but on a constructive method that is less convenient from a practical point of view: the three binomials described there are 
\begin{eqnarray*}
f&=&x_1^{\delta^{\ast}}-x_0^{{\epsilon_1}^{\ast}}x_3^{{\phi_1}^{\ast}}\\
f_1&=&x_2^{\partial p^{\mu}}-x_0^{\alpha_0}x_1^{\alpha_1}x_3^{\alpha_3}\\
f_2&=&x_2^{\partial q^{\nu}}-x_0^{\beta_0}x_1^{\beta_1}x_3^{\beta_3},
\end{eqnarray*}
(assuming, without loss of generality, that $\gcd(\delta, \epsilon_1, \epsilon_2)=1$), where  $p$ and $q$ are arbitrary distinct primes, and $\mu, \nu$, and $\alpha_i, \beta_i$ are suitable nonnegative integers to be determined so as to have $f_1,f_2\in I(C)$.  
}
\end{remark}
\begin{example}\label{example5}{\rm Consider the projective curve
$$C: C: x_0=\xi^{70},\ x_1=\xi^{66}\omega^4,\ x_2=\xi^{15}\omega^{55},\ x_3=\omega^{70}.$$
\noindent
The parametrization fulfils case (IV) in the proof of Corollary \ref{corollary2}. The corresponding matrix is
$$A=\left(\begin{array}{ccc}
x_1^{35}&x_3^2&x_0^{18}x_2^4\\&&\\
x_0^{33}&1&x_1^{20}
\end{array}\right).$$
\noindent
Hence, according to Proposition \ref{ci}, $C$ is set-theoretically defined by
\begin{eqnarray*}
f=\Delta_{12}&=&x_1^{35}-x_0^{33}x_3^2\\\\
f_1=\Delta_{13}&=&x_1^{55}-x_0^{51}x_2^4\\\\
f_2&=&x_2^{22}-x_1^{5}x_3^{17}.
\end{eqnarray*}
The method in \cite{BM}, however, for $p=2$ and $q=3$, yields binomials of higher degree, namely
\begin{eqnarray*}
f&=&x_1^{35}-x_0^{33}x_3^2\\\\
f_1&=&x_2^{64}-x_0^9x_1^5x_3^{50}\\\\
f_2&=&x_2^{162}-x_0^{30}x_1^5x_3^{127},
\end{eqnarray*}
if one takes $f_1$ and $f_2$ of the least possible degree. 
}
\end{example}
\begin{remark}{\rm We know from the work by Cattani, Curran and Dickenstein \cite{CCD} that the defining ideal of every projective monomial curve in ${\bf P}^3$ contains a complete intersection ideal of height 2. Our Proposition \ref{ci} makes this result more precise: it establishes that this ideal can be chosen in such a way that {\it it has the same radical as the determinantal ideal generated by the $2$-minors of a $2\times 3$ simple matrix with monomial entries}. In \cite{CCD} the authors consider the curve known as the {\it twisted cubic}
$$C: x_0=\xi^3,\ x_1=\xi^2\omega,\ x_2=\xi\omega^2,\ x_3=\omega^3,$$
\noindent
and find that $I(C)$ contains the ideal $(x_1^2-x_0x_2,\ x_2^3-x_0x_3^2)$.   Curve $C$ fulfils case (I) in the proof of Corollary \ref{corollary2}; hence it is associated with the matrix  
$$A=\left(\begin{array}{ccc}
x_1^2&x_0x_3&x_2\\
x_0&x_1&1
\end{array}\right).$$
\noindent
We deduce, as in Example \ref{example4}, that   
 $(x_1^2-x_0x_2,\ x_0x_3-x_1x_2)\subset I(C)$, and $(x_1^2-x_0x_2,\ x_1^3-x_0^2x_3)\subset I(C)$.} 
\end{remark}
All the above can be applied to affine monomial curves in $K^3$: since these are the affine parts of projective monomial curves in ${\bf P}^3$, it suffices to set $x_0=\xi=1$.

\end{document}